\documentclass[10pt,twocolumn,letterpaper]{article}
\usepackage{subfig}
\usepackage{cvpr}
\usepackage{times}
\usepackage{epsfig}
\usepackage{graphicx}
\usepackage{amsmath}
\usepackage{amssymb}
\usepackage{array}
\newcolumntype{P}[1]{>{\centering\arraybackslash}p{#1}}

\newcommand{\comment}[1]{}

\cvprfinalcopy 

\def\cvprPaperID{7} 

\ifcvprfinal\fi
\begin{document}

\title{A Riemannian Framework for Statistical Analysis of Topological Persistence Diagrams}

\author{Rushil Anirudh*, Vinay Venkataraman*, Karthikeyan Natesan Ramamurthy$^{+}$, Pavan Turaga*\\
*Arizona State University, $^{+}$ IBM T. J. Watson Research Center, NY \\
\and
}

\maketitle

\begin{abstract}
Topological data analysis is becoming a popular way to study high dimensional feature spaces without any contextual clues or assumptions. This paper concerns itself with one popular topological feature, which is the number of $d-$dimensional holes in the dataset, also known as the Betti$-d$ number. The persistence of the Betti numbers over various scales is encoded into a persistence diagram (PD), which indicates the birth and death times of these holes as scale varies. A common way to compare PDs is by a point-to-point matching, which is given by the $n$-Wasserstein metric. However, a big drawback of this approach is the need to solve correspondence between points before computing the distance; for $n$ points, the complexity grows according to $\mathcal{O}($n$^3)$. Instead, we propose to use an entirely new framework built on Riemannian geometry, that models PDs as 2D probability density functions that are represented in the square-root framework on a Hilbert Sphere. The resulting space is much more intuitive with closed form expressions for common operations. The distance metric is 1) correspondence-free and also 2) independent of the number of points in the dataset. The complexity of computing distance between PDs now grows according to $\mathcal{O}(K^2)$, for a $K \times K$ discretization of $[0,1]^2$. This also enables the use of existing machinery in differential geometry towards statistical analysis of PDs such as computing the mean, geodesics, classification etc. We report competitive results with the Wasserstein metric, at a much lower computational load, indicating the favorable properties of the proposed approach.
\end{abstract}

\section{Introduction}
Topological data analysis (TDA) has emerged as a useful tool to analyze underlying properties of data before any contextual modeling assumptions kick in. Generally speaking, TDA seeks to characterize the \emph{shape} of high dimensional data (viewed as a point-cloud in some metric space), by quantifying various topological constructs such as connected components, high-dimensional holes, level-sets and monotonic regions of functions defined on the data \cite{edelsbrunner2002topological}. In particular, the number of $d$-dimensional holes in a data, also known as the \textit{Betti}-$d$ number, corresponds to the rank of the $d-$dimensional homology group. A commonly used and powerful topological feature, the \textit{persistence diagram} (PD) summarizes the \textit{persistence} of the Betti numbers with respect to the \textit{scale} parameter used to analyze the data. A typical machine learning pipeline using TDA features would first estimate the PDs from the given point cloud, and define a metric on them to compare different point clouds. The Wasserstein distance measure has become ubiquitous as a metric between such PDs, as it is stable and has a well-defined metric space associated with it \cite{turner2014frechet, cohen2010lipschitz}.  However, computation of the Wasserstein distance involves finding a one-to-one map between the points in one persistence diagram to those in the other, which is a computationally expensive operation.\let\thefootnote\relax\footnotetext{This work was supported by NSF CAREER grant 1452163.}

\begin{figure*}[!tb]
\centering
\includegraphics[clip = true,trim=0mm 88mm 0mm 50mm,width = 0.99\linewidth]{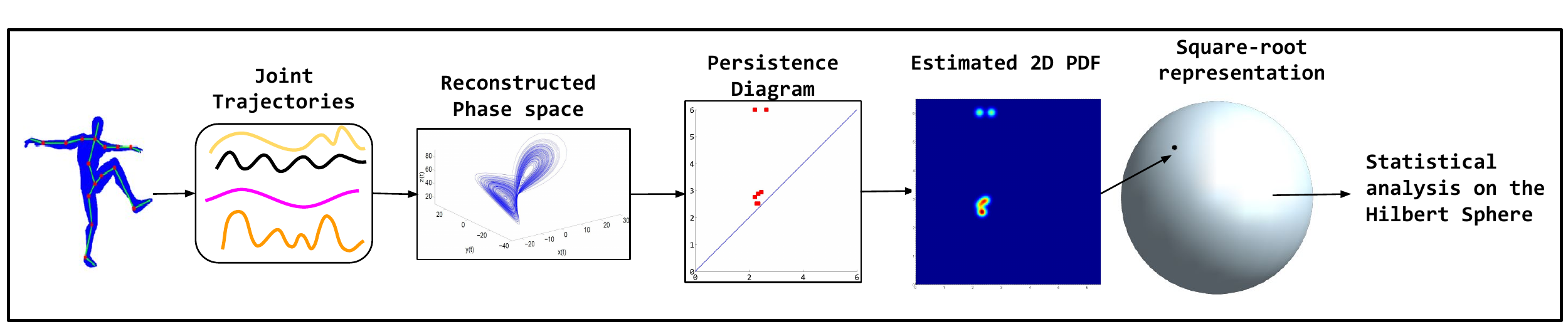}
\caption{\small{The overall sequence of operations leading to the proposed representation, illustrated for the application of activity analysis. Joint position data from motion capture systems are collected as 1D time series. The phase space is reconstructed from these time series via Takens' embedding theorem. We represent the topological properties of this phase space using the persistence diagram (PD). Next, we use kernel density estimation to represent the PD itself as a 2D probability density function (pdf). Finally, we use the square-root framework to interpret these pdfs as points on a Hilbert sphere.}}
\label{fig:pipeline}
\end{figure*}

In this paper, we propose a novel representation for persistence diagrams as points on a hypersphere, by approximating them as 2D probability density functions (pdf). We perform a square-root transform of the constructed pdf, wherein the Hilbert sphere becomes the appropriate space for defining metrics \cite{Srivastava2007}. This insight is used to construct closed form metrics -- geodesics on the Hilbert sphere -- which can be computed very efficiently, bypassing the correspondence problem entirely. The overall pipeline of operations for computing the proposed representation is given in Figure \ref{fig:pipeline}.


The biggest advantage of the proposed representation is that it completely works around the computationally expensive step of obtaining one-to-one correspondences between points in persistence diagrams, thereby making the distance computation between PDs extremely efficient. We show that approximating PDs as pdfs results in comparable performances to the popular and best performing $L_1$-Wasserstein metric, while at the same time being orders of magnitude faster. We also provide a theoretically well-grounded understanding of the geometry of the pdf representation. Additionally, the availability of closed form expressions for many important tools such as the geodesics, exponential maps etc. enables us to adapt powerful statistical tools -- such as clustering, PCA, sample mean etc. -- opening new possibilities for applications involving large datasets. To the best of our knowledge we are the first to propose this representation for persistence diagrams.

\paragraph{Contributions}
\begin{enumerate}
\item We present the first representation of persistence diagrams that are approximated as points on a Hilbert sphere resulting in closed-form expressions to compare two diagrams. This also completely avoids the correspondence problem, which is typically a computational bottleneck.
\item We demonstrate the ability of the proposed representation for statistical analyses, such as computing the mean persistence diagram, principal geodesic analysis (PGA), and classification using SVMs. 
\item We show promising results for supervised tasks such as action recognition, and assessment of movement quality in stroke rehabilitation.
\item The space of the representation -- the Hilbert sphere -- is a geometrically well-understood and intuitive space, which may further promote its use in topological machine learning applications. 
\end{enumerate}

The rest of the paper is organized as follows. Section \ref{sec:rel_work} discusses related work in more detail. Section \ref{sec:math_prelim} provides the necessary background on persistent homology, the space of persistence diagrams, and the square-root framework on the Hilbert space.  Section \ref{sec:prop_framework} provides details about the proposed framework for using the new representation of persistence diagrams in a functional Hilbert space for statistical learning tasks. Section \ref{sec:expts}  describes the experiments and results. Section \ref{sec:disc} concludes the paper.

\section{Related Work}
\label{sec:rel_work}
Persistence diagrams elegantly summarize the topological properties of point clouds, and the first algorithm for computing topological persistence was proposed in \cite{edelsbrunner2002topological}. Ever since, there has been an explosion in understanding the properties of PDs, comparing PDs, and computing statistics on them. Since a PD is a multiset of off-diagonal points along with infinite number of points on the diagonal, the cost of optimal bijective assignment, also known as the Wasserstein distance, between the individual points in the PDs is a valid distance measure. The time complexity of computing the distance between two PDs with $n$ points each is $\mathcal{O}(n^3)$ \cite{bertsekas1981new}. It has been shown in \cite{cohen2010lipschitz} that the persistence diagrams of Lipschitz functions are stable with respect to $p$-Wasserstein distance. However, the bottleneck and $p$-Wasserstein distance do not allow for easy computation of statistics. Hence, $L_p$-Wasserstein \textit{metrics} have been used to develop approaches for computing statistical summaries such as the Fr\'{e}chet mean \cite{turner2014frechet}. Computing Fr\'{e}chet mean of PDs involves iteratively assigning points from the individual diagrams to the candidate mean diagram, recomputing the mean, and repeating till convergence. Since the mean PD obtained is not guaranteed to be unique, the authors in \cite{munch2013probabilistic} proposed the  probabilistic Fr\'{e}chet mean which itself is a probability measure on the space of PDs. An investigation of the properties of the space of PDs that allow for definition of various probability measures is reported in \cite{mileyko2011probability}, and a derivation of confidence sets that allow the separation of topological signal from noise is presented in \cite{fasy2014confidence}. 

Since operations with PDs can be computationally expensive, Bubenik proposed a new topological summary - known as the persistence landscape - derived from the PD \cite{bubenik2015statistical}. It can be thought of as a collection of envelope functions on the points in PD based on the order of their importance. Persistence landscapes allow for fast computation of statistics since they lie in a vector space. However their practical utility has been limited since they provide decreasing importance to secondary and tertiary features in PDs that are usually useful in terms of discriminating between data from different classes. Recently, an approach that defines a stable multi-scale kernel between persistence diagrams has been proposed \cite{reininghaus2015stable}. The kernel is obtained by creating a surface with a Gaussian centered at each point in the PD along with a negative Gaussian centered at its reflection across the diagonal.  In addition, the authors in \cite{adams2015persistent} propose to compute a vector representation - the persistence image - by computing a kernel density estimate of the PD and integrating it locally. The kernel density estimate is weighted such that the points will have increasing weights as they move away from the diagonal. A similar weighting approach is used in \cite{kusano2016persistence} to create a persistence-weighted Gaussian kernel. 



The square-root representation for 1D probability density functions (pdfs) was proposed in \cite{Srivastava2007} and was used for shape classification. It has since been used in several applications, including activity recognition \cite{VeeraraghavanPAMI2009} - where the square-root velocity representation is used to model the space of time warping functions. This representation extends quite easily to arbitrary high-dimensional pdfs as well.

\section{Mathematical Preliminaries}
\label{sec:math_prelim}
In this section we will introduce the concepts of a) persistence homology, b) persistent diagrams, c) our proposed representation of PDs as a 2D pdf, and the square-root framework resulting in the Hilbert sphere geometry. 

\subsection{Persistent Homology}
\label{sec:persistent homology}
\begin{figure}[!htb]
\captionsetup[subfigure]{labelformat=empty}
\centering		
\subfloat[]{
\includegraphics[width=0.45\linewidth]{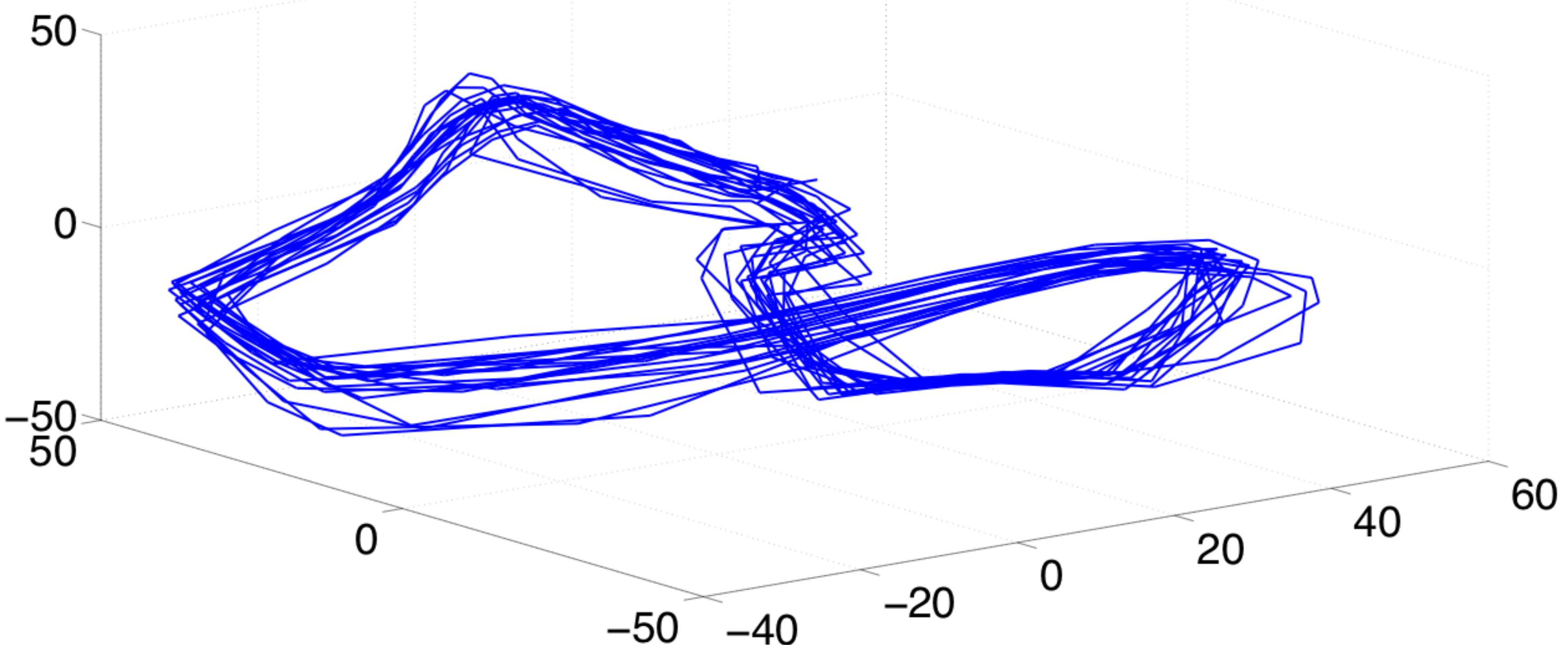}
 }
 \subfloat[]{
\includegraphics[width=0.45\linewidth]{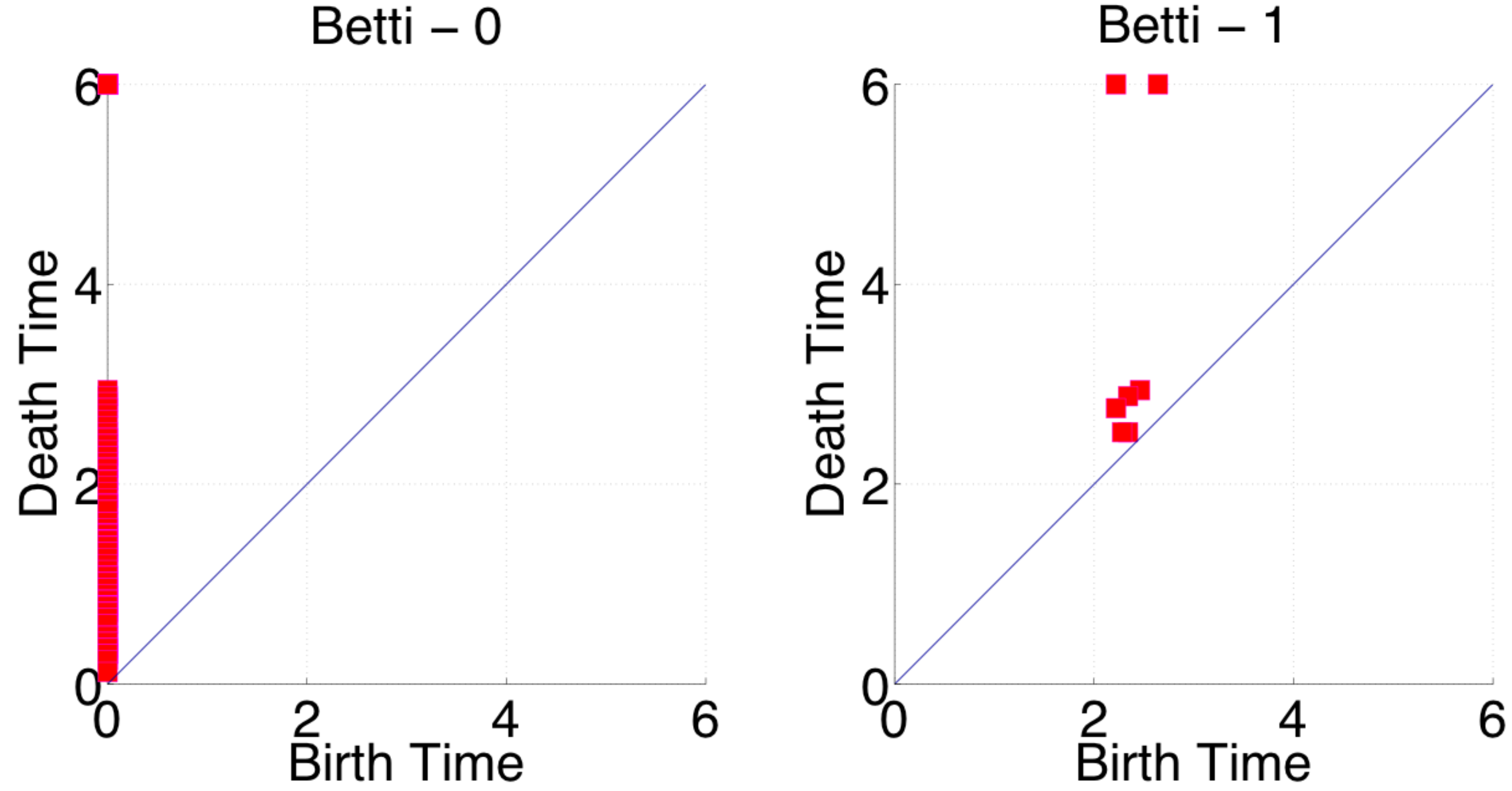}
}

\subfloat[Point cloud data]{
\includegraphics[width=0.45\linewidth]{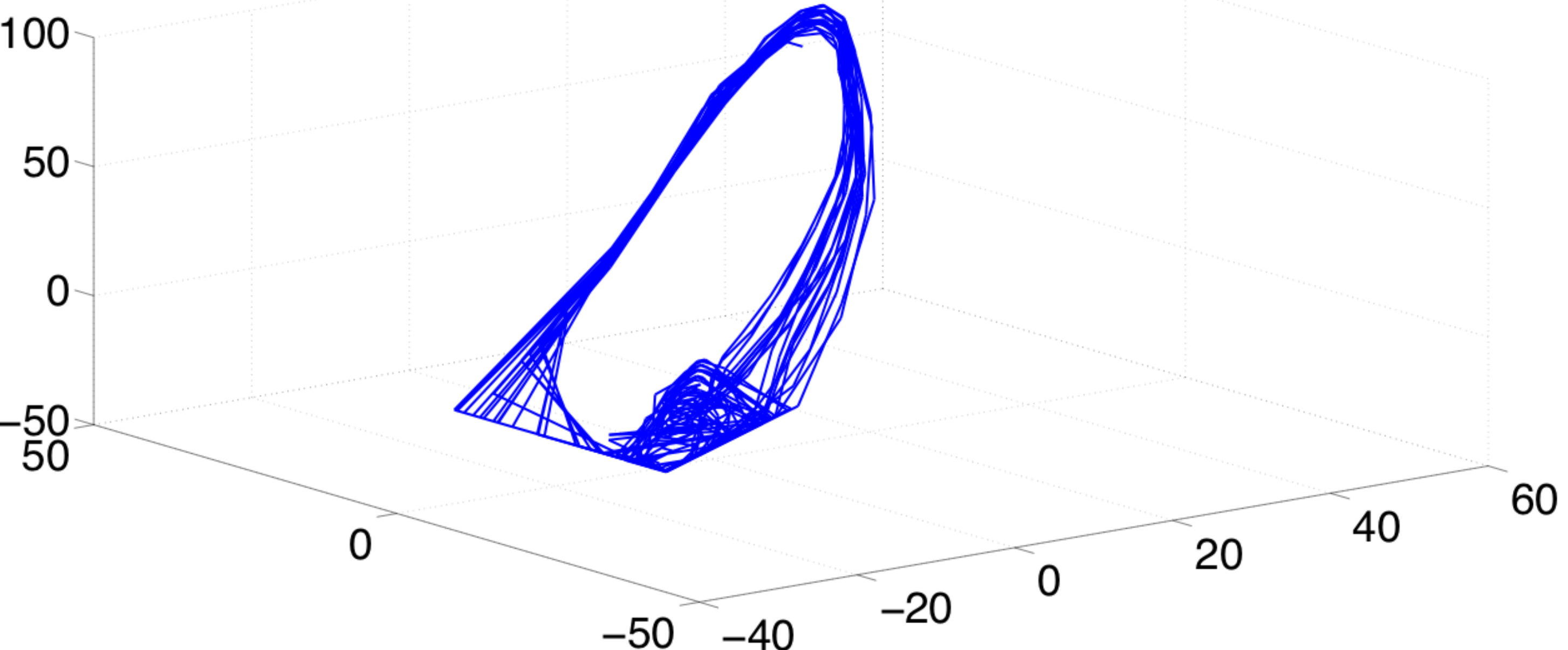}

}
\subfloat[Persistence diagrams]{
\includegraphics[width=0.45\linewidth]{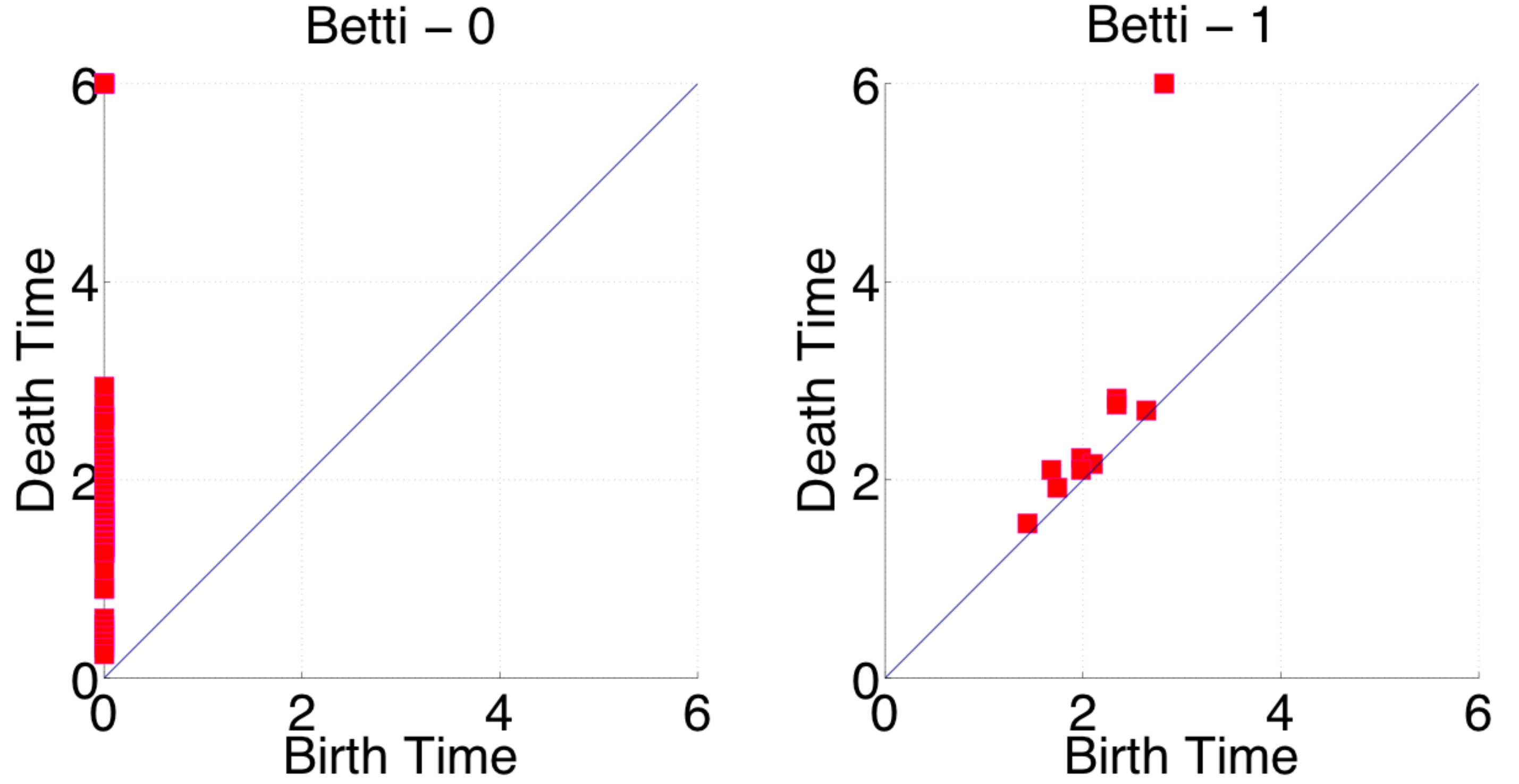}}

	\caption{\small{The above example illustrates the topological features extracted from the point cloud data with two and one one-dimensional holes. These properties are reflected well in their corresponding persistence diagrams on the right.}}
	\label{fig:PhaseSpace}
\end{figure}

The homology of point cloud $\mathbf{X}$ can be computed by first constructing a \textit{shape} out of the point cloud and estimating its homology. A simplex is imposed on every set of points in the point cloud that satisfy a neighborhood condition based on a scale parameter $\epsilon$. The collection of such simplices is known as the \textit{simplicial complex} $\mathcal{S}$. $\mathcal{S}$ can be thought of as a representation of the underlying shape of the data, and its homology can be inferred. However, homology groups obtained from $\mathcal{S}$ depend on the scale or time parameter $\epsilon$ based on which the complexes are constructed \cite{edelsbrunner2002topological}. The homological features of the simplicial complex constructed from the point cloud that are stable across scales, i.e., that are \textit{persistent}, are the ones that provide information about the underlying shape. Topological features that do not persist are considered to be noise. This information is represented using \textit{persistence diagrams} as a 2D plot of birth versus death times of each homology cycle corresponding to the homology group $\mathbb{H}_k, k = \{0, 1, \ldots\}$. The birth time is the scale at which the homological feature is born and the death time is the scale at which it ceases to exist. The homology cycle of dimension $d$ is also referred to as a $d-$dimensional hole. Therefore, the PD can be considered as an object that represents the number of holes in terms of the range of scales at which they appear. Typically, PDs of the point cloud data are obtained using \textit{filtration} of simplicial complexes. A well-known filtration is the Vietoris-Rips (VR) filtration, where a simplex is induced between a group of points whenever the distance between each pair of points is less than the given scale $\epsilon$ \cite{zomorodian2010fast}. An example of point clouds and their corresponding persistence diagrams for homology groups 0 and 1 are provided in Figure \ref{fig:PhaseSpace}.

\subsection{The Space of Persistence Diagrams}
\label{sec:pers_diag_space}
Every PD is a multiset of 2D points, where each point denotes the birth and death time of a homology group. Furthermore, the diagonal is assumed to contain an infinite number of points with the same birth and death times. For any two PDs $X$ and $Y$, the distance between the diagrams is usually quantified using the bottleneck distance or the $L_q$-Wasserstein metric \cite{kerbergeometry}. In this paper, we consider only the $L_2$- and $L_1$-Wasserstein ($d_{L_2}$ and $d_{L_1}$) metrics given respectively as,
\begin{equation}
	d_{L_2}(X,Y) = \left( \inf_{\eta:X \to Y}  \sum_{x \in X} ||x-\eta(x)||_2^2 \right)^{\frac{1}{2}},
\end{equation} and
\begin{equation}
	d_{L_1}(X,Y) = \inf_{\eta:X \to Y}  \sum_{x \in X} ||x-\eta(x)||_1.
\end{equation} Since each diagram contains an infinite number of points in the diagonal, this distance is computed by pairing each point in one diagram uniquely to another non-diagonal or diagonal point in the other diagram, and then computing the distance. This correspondence can be obtained via the Hungarian algorithm or its more efficient variants \cite{kerbergeometry}. 

The space of PDs with respect to the $L_2$-Wasserstein metric is given by
\begin{equation}
\mathcal{D}_{L_2} = \{X | d_{L_2}(X,\varnothing) < \infty\}.
\end{equation} Turner {\em et. al.} \cite{turner2014frechet} show this is a non-negatively curved Alexandrov space. Furthermore, the diagram on the geodesic between the PDs $X$ and $Y$ in this space is given as
\begin{equation}
\gamma(s) = (1-s) x + s \phi (x),
\end{equation} where $x$ is a point in the diagram $X$, $\phi(x)$ is a corresponding point in the diagram $Y$ and $s \in [0,1]$ parametrizes the geodesic. Clearly, the points in the diagram on the geodesic can be simply obtained by linearly interpolating between the corresponding points on the candidate points $X$ and $Y$. Furthermore, the Riemannian mean between two PDs is easily computed as $\gamma(0.5)$. 

\subsection{Square-root Framework and the Hilbert Sphere}
\label{sec:sqrt_fmwk}
We treat the points in a persistence diagram as samples from an underlying probability distribution function. This representation is inspired from recent work dealing with defining Mercer kernels on PDs \cite{reininghaus2015stable}, where the feature embedding obtained from the Mercer kernel closely resembles a kernel-density estimator from the given points, with an additional reflection term about the diagonal to account for boundary effects while solving a heat-equation. The fact that the feature embedding resembles a kernel density estimate is not further exploited in \cite{reininghaus2015stable}.

In our work, we more directly exploit this pdf interpretation of PDs, further endowing it with a square-root form \cite{Srivastava2007} -- and then utilizing the Hilbert spherical geometry that results. 

Without loss of generality, we will assume that the supports for each 2D pdf is $[0,1]^2$. The space of pdfs that we will consider is
\begin{multline}
\mathcal{P} = \{p:[0,1]\times [0,1]\rightarrow \mathbb{R} ~~ \forall x,y| p(x,y)\geq 0, \\
\mbox{ and  } \int_0^1\int_0^1p(x,y)dxdy = 1\}
\end{multline} It is well-known that $\mathcal{P}$ is not a vector space \cite{Srivastava2007}. Instead, it is a Riemannian manifold with the Fisher-Rao metric as the natural Riemannian metric. Geodesics under the Fisher-Rao metric are quite difficult to compute. Instead, we adopt a square-root form proposed by Srivastava {\em et. al.} \cite{Srivastava2007} which simplifies further analysis enormously. In other words we consider the space,
\begin{multline}
\Psi = \{\psi:[0,1]\times[0,1] \rightarrow \mathbb{R}| \psi\geq 0,\\
\int_0^1\int_0^1 \psi^2(x,y)dxdy = 1\}.
\end{multline} For any two tangent vectors $v_1,v_2 \in T_\psi(\Psi)$, the Fisher-Rao metric is given by,
\begin{equation}
\langle v_1,v_2\rangle = \int_0^1\int_0^1v_1(x,y)v_2(x,y)~dx dy.
\end{equation}

The above two pieces of information imply that the square-root form $\psi = \sqrt{p}$ results in the space becoming a unit Hilbert-sphere, endowed with the usual inner-product metric. Geodesics on the unit-Hilbert sphere under the above Riemannian metric are known in closed form. In fact, the differential geometry of the Hilbert sphere results in closed form expressions for computing geodesics, exponential and inverse exponential maps \cite{Srivastava2007}. Further, the square-root form with the above metric has additional favorable properties such as invariance to domain re-parameterization. 

Given two points $\psi_1, \psi_2$ the geodesic distance between them is given by
\begin{equation}
\label{eqn:hilb_dist}
d_H(\psi_1,\psi_2) = \cos^{-1}(\langle \psi_1,\psi_2\rangle),
\end{equation}
where $\langle \psi_1,\psi_2\rangle_\psi = \int_0^1 \psi_1(t)\psi_2(t) dt$. The exponential map is defined as,
 \begin{equation}
\mathrm{ exp}_{\psi_i}(\upsilon) = \mathrm{cos}(||{\upsilon}||_{\psi_i})\psi_i + \mathrm{sin}(||\upsilon||_{\psi_i})\frac{\upsilon}{||\upsilon||_{\psi_i}},
 \end{equation} where $\upsilon \in T_{\psi_i}( \Psi)$ is a tangent vector at $\psi_i\mbox{ and }||\upsilon||_{\psi_i}$ = $\sqrt{\left<\upsilon,\upsilon\right>_{\psi_i}}$. 
The logarithmic map from $\psi_i$ to  $\psi_j$ is given by:
\begin{equation}
\exp^{-1}_{\psi_i}(\psi_j) = \frac{u}{\sqrt{\langle u,u\rangle}} \mathrm{cos}^{-1}\left< \psi_i, \psi_j\right>,
\end{equation}
with $u = \psi_i - \left< \psi_i, \psi_j\right>\psi_j.$ The geodesic on the sphere is given in closed form as
\begin{equation}
\label{eqn:geod1}
\pi(s) = \frac{(1-s)\psi_1+s \psi_2}{s^2+(1-s)^2+2s(1-s)(\langle \psi_1, \psi_2\rangle)},
\end{equation} or equivalently as
\begin{equation}
\label{eqn:geod2}
\pi(s) = \cos(s||v||) \psi + \sin(s||v||)\frac{v}{||v||}.
\end{equation}

\begin{figure*}[!htb]
\centering
\includegraphics[clip = true,trim=0mm 5mm 0mm 0mm,width = 0.95\linewidth]{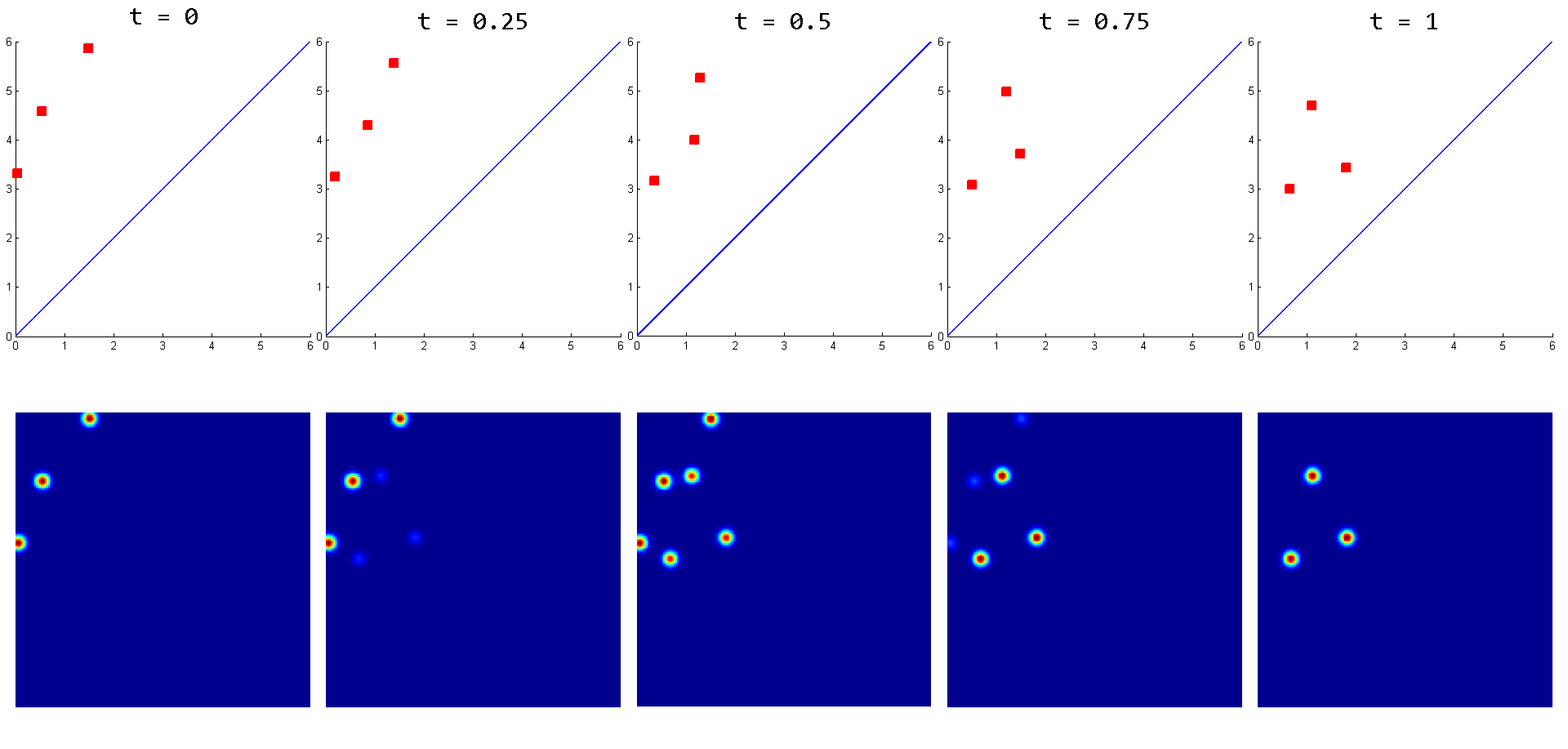}
\caption{\small{Comparing geodesics in the Alexandrov space (top) and the proposed Hilbert sphere (bottom), for a simple persistence diagram containing just $3$ points. The Alexandrov space requires computing correspondence between points, and the geodesic involves linearly interpolating between corresponding points. The persistence pdf avoids correspondence estimation, hence the geodesics correspond to new local modes appearing and gradually increasing in intensity as the original modes decrease in intensity.}}
\label{fig:sampling}
\end{figure*}

A comparison of sampling of the geodesics in the Alexandrov space induced by the $L_2$-Wasserstein metric and the proposed Hilbert sphere are provided in Figure \ref{fig:sampling}.

\section{Algorithmic Details}
\label{sec:prop_framework}
Our proposed framework consists of reconstructing the phase space of the time series data, computing the PDs of the reconstructed phase space, transforming each PD to a 2D pdf by kernel density estimation, using the 2D pdfs in the square-root framework to estimate distances or obtain statistical summaries. The distances computed and the features obtained will be used in inference tasks. Additional details for some of the above steps are given below.

\paragraph{Phase-space Reconstruction from Activity Data}
In dynamical system understanding, we aim to estimate the underlying states, but we measure functions -- usually projections -- of the original state-space. e.g. while human movement is influenced by many joints, ligaments, muscles of the body, we might measure the location of only a few joints. One way to deal with this issue is to reconstructing the `phase space' by the method of delays used commonly in non-linear dynamical analysis \cite{Takens}. Reconstructing the phase-space in this manner only preserves important topological properties of the original dynamical system, and does not recover the true state-space. For a discrete dynamical system with a multidimensional phase-space, time-delay vectors (or embedding vectors) are obtained by the concatenation of delayed versions of the data points,
\begin{equation}
	\mathbf{x}_t = [x(t),x(t+\tau),\cdots,x(t+(m-1)\tau)]^T. 
	\label{EmVec}
\end{equation} where $m$ is the embedding dimension, $\tau$ is the embedding delay. For a sufficiently large $m$, the important topological properties of the unknown multidimensional system are reproduced in the reconstructed phase-space \cite{abarbanel1996analysis}. The collection of $n$ time-delay vectors is the point cloud data under further consideration, and this is represented as $\mathbf{X} = [\mathbf{x}_t]_{t=0}^n$.

\paragraph{Estimating the PDs}
After estimating the point cloud in the reconstructed phase space, we will construct a simplicial complex $\mathcal{S}$ using the point cloud data $\mathbf{X}$ to compute the persistence diagrams for the Betti numbers using the VR filtration. However, this approach considers only spatial adjacency between the points and ignores the temporal adjacency. We improve upon the existing VR filtration approach by explicitly creating temporal links between $\mathbf{x}_{t-1}$, $\mathbf{x}_{t}$, and $\mathbf{x}_{t+1}$ in the one-skeleton of $\mathcal{S}$, thereby creating a metric space which encodes adjacency in both space and time \cite{venkataraman2016persistent}. The persistence diagrams for homology groups of dimensions $0$ and $1$ are then estimated.

\paragraph{Square-root Framework for Distance Estimation between PDs}
Since each PD is a multiset of points in the 2D space, we start by constructing a 2D pdf from each of them using kernel density estimation using a Gaussian kernel of zero mean and variance $\sigma^2$. For each PD, we compute the square-root representation $\psi$ using the approach provided in Section \ref{sec:sqrt_fmwk}, and the distance between two PDs can be computed using (\ref{eqn:hilb_dist}).

\paragraph{Dimensionality Reduction with Principal Geodesic Analysis (PGA)}
We are able to extract principal components using PGA \cite{Fletcher2004}, adapted to our Hilbert sphere -- which essentially performs classical PCA on the tangent space of the mean-point on the sphere. Given a dataset of persistence diagrams in their square-root form $\{\psi_1,\psi_2,\ldots, \psi_N\}$, we first compute the Riemannian center of mass \cite{GroveK72}. We use a simple approximation using an extrinsic algorithm that computes first the Euclidean mean and maps it on to the closest point on the sphere. Next, we represent each $\psi_i$ by its tangent vector from the mean. We then compute the principal components on the tangent space and represent a persistence diagram as a low-dimensional vector.

\section{Experiments}
We perform experiments on two datasets for human action analysis. First we perform action recognition on the MoCap dataset \cite{ali2007chaotic}, next we show the use of the framework in quality assessment of movements in stroke rehabilitation \cite{Chen2011Stroke}. We will describe the datasets next, followed by the evaluation settings and parameters used. In all our experiments, we performed a discretization of the 2D pdf into a $K\times K$ grid. The choice of $K$ indirectly affects classification performance as expected, i.e. -- a smaller value of $K$ results in a reduced ability to identify between nearby points in the PD due to lower resolution. On the other hand, a larger value of $K$ improves resolution, but also increases computational requirements. Typical values of $K$ used in experiments range from $50$ to $100$ at most, whereas typical values of $n$ -- the number of points in a PD -- ranges from $20$ to $50$.

\label{sec:expts}
\subsection{Motion Capture Data}
We evaluate the performance of the proposed framework using $3$-dimensional motion capture sequences of body joints \cite{ali2007chaotic}. The dataset is a collection of five actions: \textit{dance, jump, run, sit} and \textit{walk} with $31, 14, 30, 35$ and $48$ instances respectively. 
\begin{table}[!htb]
\centering
		\begin{tabular}{|p{1in}|c|p{1.0in}|}
			\hline
			\textbf{Method} &\small{\textbf{Accuracy (\%)}} & \textbf{Time (sec)}	\\ \hline \hline
			Chaos \cite{ali2007chaotic} & 52.44 & - \\ 
			VR-Complex \cite{zomorodian2010fast} & 93.68 & -   \\ 
			DT2 \cite{vinay_PAMI} & 93.92 & - 	\\ 
            T-VR Complex (L1) \cite{venkataraman2016persistent} & 96.48 & $(1.2 \pm 1.23) \times 10^3 $  \\ \hline
            \textbf{Proposed} - 1NN & {89.87} & $\mathbf{(0.059 \pm 0.044)}$ \\
			\textbf{Proposed} - PGA +SVM & {91.68}& -  \\ \hline
		\end{tabular}
        \vspace{10pt}
\caption{\small{Comparison of classification rates for different methods using nearest neighbor classifier on the motion capture dataset. It is observed that on an average, the proposed metric is $10^5$ times faster than the traditional Wasserstein metric, while achieving a comparable recognition accuracy.}}
		\label{tab:results1}
        
\end{table}
We generate $100$ random splits having $5$ testing examples from each action class and use an SVM classifier on the vector features computed with PGA, we get a performance of $91.68\%$. The mean recognition rates for the different methods are given in Table \ref{tab:results1}. We also compare with a 1-nearest neighbor classifier computed on the Hilbert sphere, which gives a performance of $89.87\%$. This is clearly competitive, when compared to approaches that use the $L_1$-Wasserstein metric as shown in table \ref{tab:results1}. Clearly, the proposed metric for persistence diagrams is able to capture most of the important information, with the added benefit being free from expensive correspondence estimation. To compute the times taken using the Wasserstein metric and the proposed metric, we average across 3 samples from 5 action classes, across $57$ (19 joint trajectories along 3 axes) trajectories -- a total of $855$ persistence diagrams. We computed a pairwise distance matrix for all these PDs, and computed the average time taken in Matlab 2014 on a standard Intel i7 PC.

\begin{figure}[!htb]
\centering
\includegraphics[clip = true,trim=150mm 0mm 140mm 0mm,width = 0.95\linewidth]{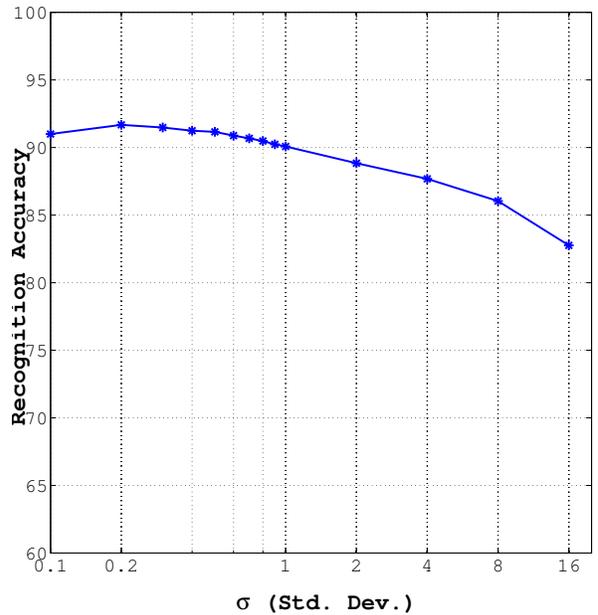}
\caption{\small{Recognition accuracy vs $\sigma$, the standard deviation of the 2D Gaussians used in the kernel density estimates.}}
\label{fig:sigma}
\vspace{-0.1in}
\end{figure}

In Figure \ref{fig:sigma}, we compare the recognition accuracy with the choice of $\sigma$, the standard deviation used during kernel density estimation. The accuracy is generally higher for a small $\sigma$ and drops as the $\sigma$ increases. We note that a similar trend is also reported in \cite{reininghaus2015stable}.

\begin{figure*}[!htb]
\centering
\subfloat[\footnotesize{Unimpaired Subjects}]{
  \includegraphics[clip = true,trim=5mm 5mm 10mm 10mm,width = 0.45\linewidth]{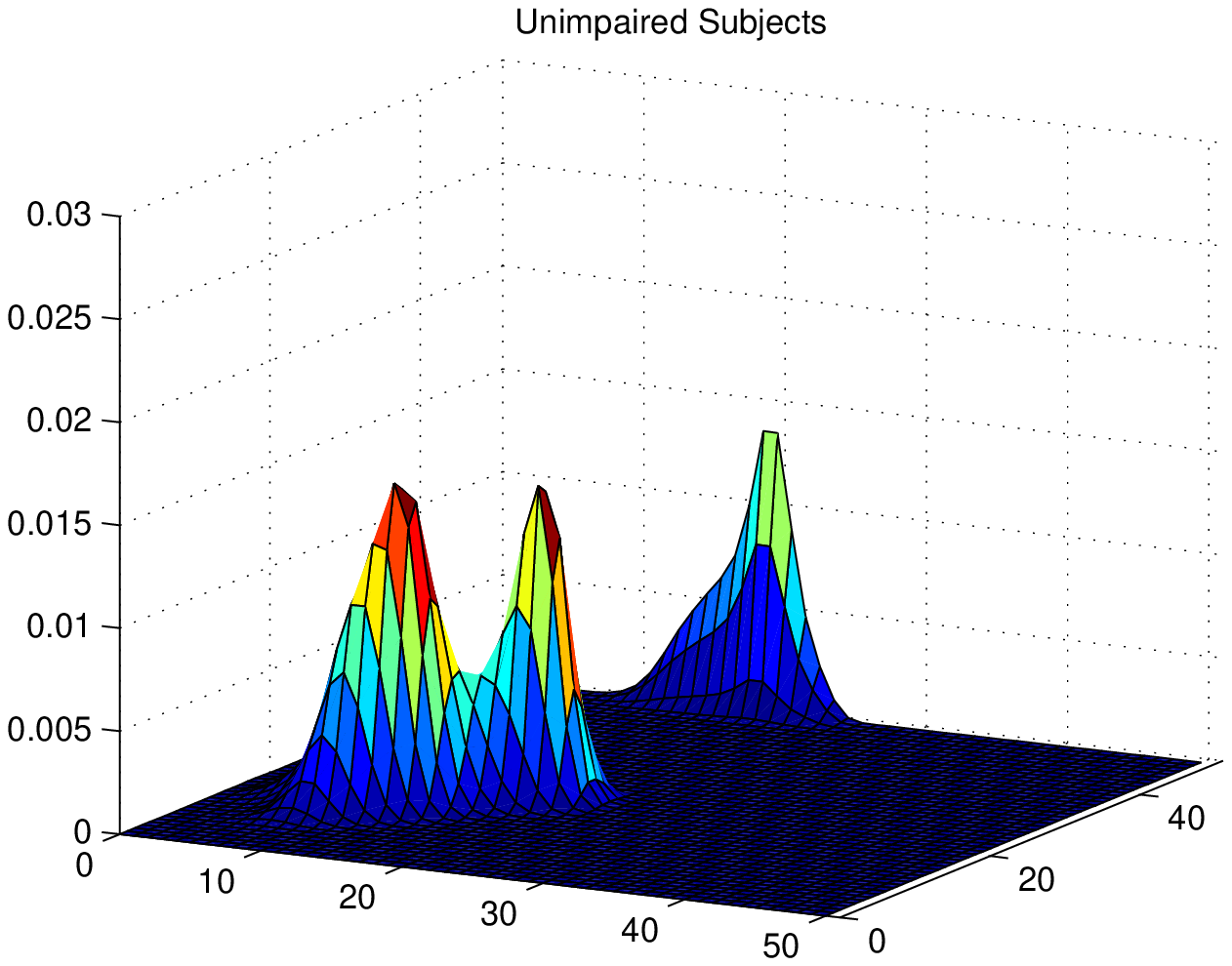}
  \label{fig:unimp_HT}}
\subfloat[\footnotesize{Impaired Subjects}]{
  \includegraphics[clip = true,trim=5mm 5mm 10mm 10mm,width = 0.45\linewidth]{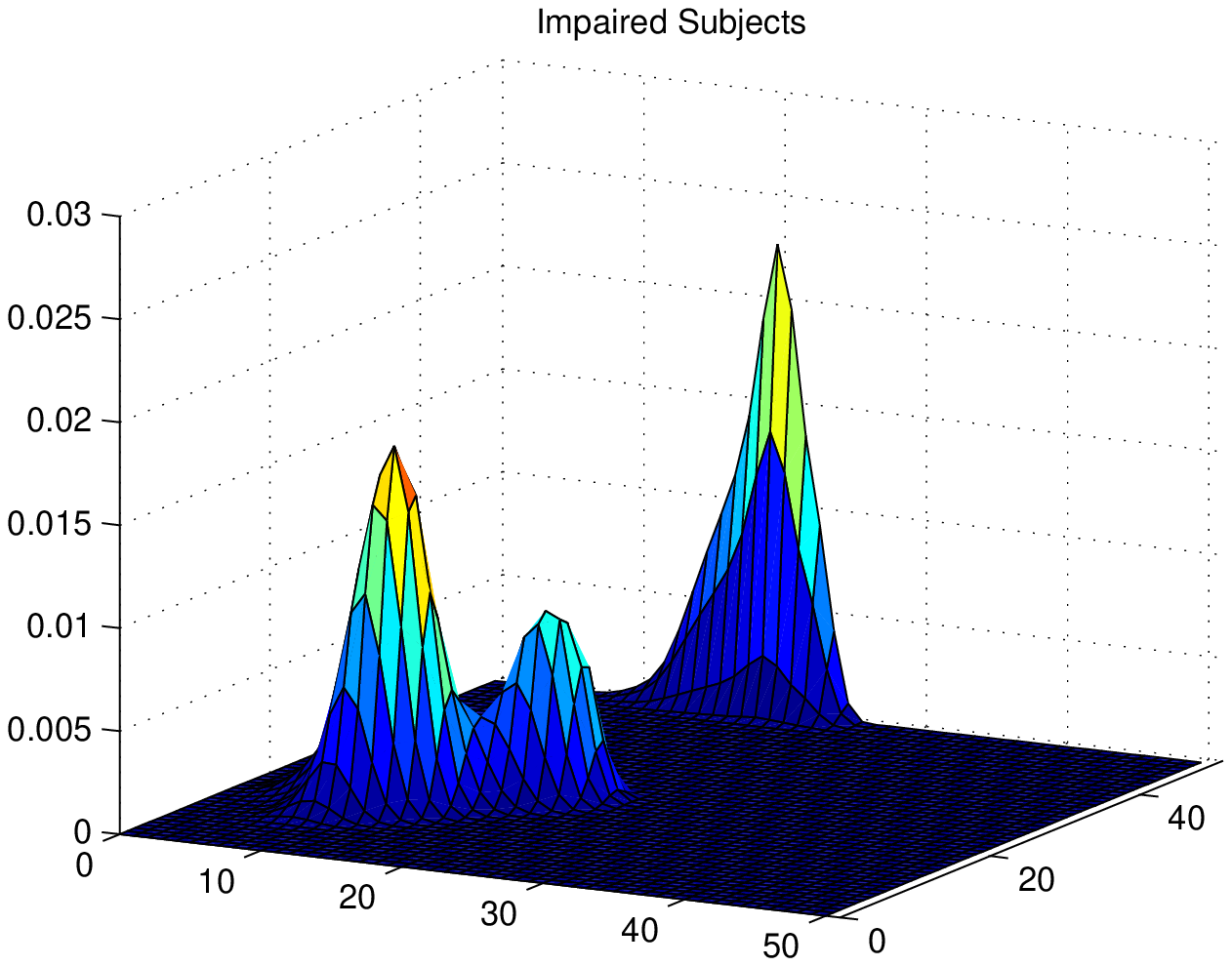}
  \label{fig:imp_HT}}
  \caption{\small{The mean persistence pdfs visualized as heat maps for unimpaired subjects (left) and stroke survivors (right), for the reach-and-grasp action, are shown. These means are computed as the extrinsic means (we choose extrinsic mean here for simplicity) on the Hilbert sphere using the proposed representation. The locations of the peaks are exactly the same, since the subjects perform the same general movement, however the intensity differs significantly indicating that the \emph{quality} of movement is captured using topological tools well}.}
\label{fig:stroke_mean}
\end{figure*}

\subsection{Stroke Rehabilitation Dataset}
Our aim in this experiment is to quantitatively assess the quality of movement performed by the impaired subjects during repetitive task therapy. The experimental data was collected using a heavy marker-based system ($14$ markers on the right hand, arm and torso) in a hospital setting from $15$ impaired subjects performing reach and grasp movements to a defined target. The stroke survivors were also evaluated by the Wolf Motor Function Test (WMFT) \cite{wolf2001assessing} on the day of recording, which evaluates a subject's functional ability on a scale of $1-5$ (with $5$ being least impaired and $1$ being most impaired) based on predefined functional tasks. In our experiments, we only use the data corresponding to the single marker on the wrist from the heavy marker-based hospital system, which allows us to evaluate the framework in a home-based rehabilitation system setting. We utilize the WMFT scores as an approximate high-level quantitative measure for movement quality (ground-truth labels) of impaired subjects performing the reach and grasp movements.

\begin{figure}[!htb]
\centering
\includegraphics[clip = true,trim=130mm 0mm 140mm 10mm,width = 0.85\linewidth]{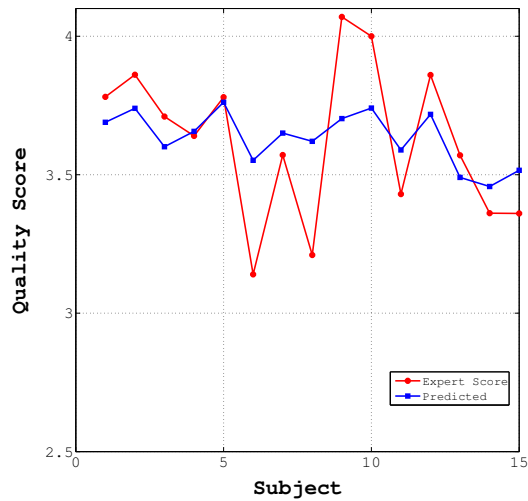}
\caption{\small{Predicting movement quality scores for reach-and-grasp in stroke survivors using topological features. We obtain a \textbf{0.8} correlation score with a p-value of $5.46 \times 10^{-4}$}}
\label{fig:stroke}
\vspace{-0.2in}
\end{figure}

We compute explicit features as described earlier using PGA. This allows us to represent each movement by a low-dimensional feature vector. We use the WMFT scores and split the dataset into train and test to perform regression. Since the dataset is small, we use a leave-one-out approach, where we train on all but one sample, which is used for testing, and repeat this over all the samples in the dataset. The correlation performance with the WMFT score is shown in Table \ref{tab:strkcomptab}, and it can be seen that we are comparable to the state of the art, and much better than traditional features. The predicted scores are shown in Figure \ref{fig:stroke} compared to the original scores. 

The dynamical features proposed in \cite{vinay_PAMI} and \cite{ali2007chaotic}, depend on describing the phase space for each movement. We compute the topological features, which are expected to be much weaker than other characteristics such as shape. However, we are still able to capture subtle information regarding movement quality. This is illustrated in Figure \ref{fig:stroke_mean}, where we see the 3D peaks associated with the average of persistence diagrams across subjects impaired with stroke and those who are unimpaired. It is clearly seen that since they are performing the same kind of movements, the peaks occur at the same location. However, the intensity is significantly different across these diagrams, perhaps indicating information regarding movement quality. 
\begin{table}[!htb]
\centering
\begin{tabular}{|c|P{0.8in}|}
\hline
\textbf{Method} & \textbf{Correlation with WMFT} 	\\ \hline
 KIM (14 markers) \cite{chen2011computational} & 0.85 		\\ 
Lyapunov exponent (1 marker) \cite{lyaprosen} & 0.60	\\ \hline 
\textbf{Proposed method (1 marker)} & \textbf{0.80}		\\ \hline
\end{tabular}
\vspace{10pt}
\caption{\small{Comparison of correlation coefficient for different methods using leave-one-subject-out cross-validation scheme and SVM regressor on the stroke rehabilitation dataset. Even with just a single marker, We obtain comparable results to a clinical 14-marker system. }}
\label{tab:strkcomptab}
\vspace{-15pt}
\end{table}

\section{Discussion and Conclusion}
\label{sec:disc}
Based on the theory and experiments presented so far, it is instructive to compare the space of PDs with respect to the different distance metrics. We will consider only the $L_2$-Wasserstein metric ($d_{L_2}$) and the proposed Hilbert sphere metric $d_{H}$. The interpretation of PDs with respect to these two metrics is very different. $d_{L_2}(X,Y)$ is the Earth Mover's Distance between the PDs $X$ and $Y$ considered as scaled discrete probability mass functions (pmfs). However, $d_{H}(\psi_X,\psi_Y)$ is the geodesic distance between the square root of kernel density estimates of $X$ and $Y$. Furthermore, the ``length'' of the persistence diagrams induced by these metrics are very different. For the Wasserstein metric, this is given by $d_{L_2}(X,\varnothing)$ which can be arbitrarily high, whereas with the proposed metric it is constrained as the persistence diagrams live on a unit sphere. The most important distinction arises when the pairwise distances between all the points in the two PDs are  sufficiently high. When the 2-dimensional pdfs obtained from kernel density estimates of the two PDs do not overlap anywhere, $d_{H}(X,Y) \rightarrow 1$. This implies that if the variance of the kernel is sufficiently small, two PDs with non-overlapping points will always have $d_H$ close to $1$. This problem can be alleviated by using kernels with multiple scales for smoothing the PDs to obtain and combining the distances obtained at each scale.

The comparison between the geodesics in the persistence space is also illuminating (Figure \ref{fig:sampling}). Sampling of the geodesic in the Alexandrov space shows that the points in one diagram move towards the corresponding points in the other, as we move along the geodesic. Whereas, a similar sampling in the Hilbert sphere shows that the PDs in the middle of the geodesic contain the modes from pdfs corresponding to both the candidate PDs. This results in markedly different interpretations of the means computed from $d_{L_2}$ and $d_H$ even for the case of two persistence diagrams.

Several directions of future work emerge from this endeavor. On the theoretical side a more formal understanding of the variety of metrics and their relationship to the proposed one can be considered. On the applied side, the favorable computational load reduction should open new applications of TDA for massive datasets. 

{\small
\bibliographystyle{ieee}
\bibliography{ref1,ref2,strings,egbib}

\begin{thebibliography}{10}\itemsep=-1pt

\bibitem{abarbanel1996analysis}
H.~D. Abarbanel.
\newblock {\em Analysis of observed chaotic data}.
\newblock New York: Springer-Verlag, 1996.

\bibitem{adams2015persistent}
H.~Adams, S.~Chepushtanova, T.~Emerson, E.~Hanson, M.~Kirby, F.~Motta,
  R.~Neville, C.~Peterson, P.~Shipman, and L.~Ziegelmeier.
\newblock Persistent images: A stable vector representation of persistent
  homology.
\newblock {\em ArXiv preprint arXiv:1507.06217}, 2015.

\bibitem{ali2007chaotic}
S.~Ali, A.~Basharat, and M.~Shah.
\newblock Chaotic invariants for human action recognition.
\newblock In {\em IEEE International Conference on Computer Vision}, pages
  1--8, Oct. 2007.

\bibitem{bertsekas1981new}
D.~P. Bertsekas.
\newblock A new algorithm for the assignment problem.
\newblock {\em Mathematical Programming}, 21(1):152--171, 1981.

\bibitem{bubenik2015statistical}
P.~Bubenik.
\newblock Statistical topological data analysis using persistence landscapes.
\newblock {\em The Journal of Machine Learning Research}, 16(1):77--102, 2015.

\bibitem{chen2011computational}
Y.~Chen, M.~Duff, N.~Lehrer, H.~Sundaram, J.~He, S.~L. Wolf, and T.~Rikakis.
\newblock A computational framework for quantitative evaluation of movement
  during rehabilitation.
\newblock In {\em AIP Conference Proceedings-American Institute of Physics},
  volume 1371, pages 317--326, 2011.

\bibitem{Chen2011Stroke}
Y.~Chen, M.~Duff, N.~Lehrer, H.~Sundaram, J.~He, S.~L. Wolf, T.~Rikakis, T.~D.
  Pham, X.~Zhou, H.~Tanaka, et~al.
\newblock A computational framework for quantitative evaluation of movement
  during rehabilitation.
\newblock In {\em AIP Conference Proceedings-American Institute of Physics},
  volume 1371, page 317, 2011.

\bibitem{cohen2010lipschitz}
D.~Cohen-Steiner, H.~Edelsbrunner, J.~Harer, and Y.~Mileyko.
\newblock Lipschitz functions have l p-stable persistence.
\newblock {\em Foundations of Computational Mathematics}, 10(2):127--139, 2010.

\bibitem{edelsbrunner2002topological}
H.~Edelsbrunner, D.~Letscher, and A.~Zomorodian.
\newblock Topological persistence and simplification.
\newblock {\em Discrete and Computational Geometry}, 28(4):511--533, 2002.

\bibitem{fasy2014confidence}
B.~T. Fasy, F.~Lecci, A.~Rinaldo, L.~Wasserman, S.~Balakrishnan, A.~Singh,
  et~al.
\newblock Confidence sets for persistence diagrams.
\newblock {\em The Annals of Statistics}, 42(6):2301--2339, 2014.

\bibitem{Fletcher2004}
P.~T. Fletcher, C.~Lu, S.~M. Pizer, and S.~C. Joshi.
\newblock Principal geodesic analysis for the study of nonlinear statistics of
  shape.
\newblock {\em IEEE Transactions on Medical Imaging}, 23(8):995--1005, August
  2004.

\bibitem{GroveK72}
K.~Grove and H.~Karcher.
\newblock How to conjugate {C$^1$}-close group actions.
\newblock {\em Math.Z}, 132:11--20, 1973.

\bibitem{kerbergeometry}
M.~Kerber, D.~Morozov, and A.~Nigmetov.
\newblock Geometry helps to compare persistence diagrams.
\newblock In {\em Proceedings of the Eighteenth Workshop on Algorithm
  Engineering and Experiments (ALENEX)}, pages 103--112, 2016.

\bibitem{kusano2016persistence}
G.~Kusano, K.~Fukumizu, and Y.~Hiraoka.
\newblock Persistence weighted gaussian kernel for topological data analysis.
\newblock {\em arXiv preprint arXiv:1601.01741}, 2016.

\bibitem{mileyko2011probability}
Y.~Mileyko, S.~Mukherjee, and J.~Harer.
\newblock Probability measures on the space of persistence diagrams.
\newblock {\em Inverse Problems}, 27(12):124007, 2011.

\bibitem{munch2013probabilistic}
E.~Munch, P.~Bendich, K.~Turner, S.~Mukherjee, J.~Mattingly, and J.~Harer.
\newblock Probabilistic {Fr{\'e}chet} means and statistics on vineyards.
\newblock {\em ArXiv preprint arXiv:1307.6530}, 2013.

\bibitem{reininghaus2015stable}
J.~Reininghaus, S.~Huber, U.~Bauer, and R.~Kwitt.
\newblock A stable multi-scale kernel for topological machine learning.
\newblock In {\em Proceedings of the IEEE Conference on Computer Vision and
  Pattern Recognition}, pages 4741--4748, 2015.

\bibitem{lyaprosen}
M.~Rosenstein, J.~Collins, and C.~De~Luca.
\newblock A practical method for calculating largest lyapunov exponents from
  small data sets.
\newblock {\em Physica D: Nonlinear Phenomena}, 65(1):117--134, 1993.

\bibitem{Srivastava2007}
A.~Srivastava, I.~Jermyn, and S.~Joshi.
\newblock Riemannian analysis of probability density functions with
  applications in vision.
\newblock In {\em IEEE Conference on Computer Vision and Pattern Recognition},
  pages 1--8, 2007.

\bibitem{Takens}
F.~Takens.
\newblock Detecting strange attractors in turbulence.
\newblock {\em Dynamical Systems and Turbulence}, 898:366--381, 1981.

\bibitem{turner2014frechet}
K.~Turner, Y.~Mileyko, S.~Mukherjee, and J.~Harer.
\newblock Fr{\'e}chet means for distributions of persistence diagrams.
\newblock {\em Discrete \& Computational Geometry}, 52(1):44--70, 2014.

\bibitem{VeeraraghavanPAMI2009}
A.~Veeraraghavan, A.~Srivastava, A.~K. Roy-Chowdhury, and R.~Chellappa.
\newblock Rate-invariant recognition of humans and their activities.
\newblock {\em Image Processing, IEEE Transactions on}, 18(6):1326--1339, 2009.

\bibitem{venkataraman2016persistent}
V.~Venkataraman, K.~N. Ramamurthy, and P.~Turaga.
\newblock Persistent homology of attractors for action recognition.
\newblock {\em arXiv preprint arXiv:1603.05310}, 2016.

\bibitem{vinay_PAMI}
V.~Venkataraman and P.~Turaga.
\newblock Shape distributions of nonlinear dynamical systems for video-based
  inference.
\newblock {\em IEEE Transactions on Pattern Analysis and Machine Intelligence},
  2016 (accepted) ; arXiv:1601.07471.

\bibitem{wolf2001assessing}
S.~L. Wolf, P.~A. Catlin, M.~Ellis, A.~L. Archer, B.~Morgan, and A.~Piacentino.
\newblock Assessing wolf motor function test as outcome measure for research in
  patients after stroke.
\newblock {\em Stroke}, 32(7):1635--1639, 2001.

\bibitem{zomorodian2010fast}
A.~Zomorodian.
\newblock Fast construction of the {Vietoris-Rips} complex.
\newblock {\em Computers \& Graphics}, 34(3):263--271, 2010.

\end{thebibliography}
}
\end{document}